\documentclass{article}

  \usepackage[colorlinks=true,linkcolor=blue,urlcolor=black,citecolor=blue,bookmarksopen=true]{hyperref}

  \usepackage[bottom]{footmisc}

  \usepackage{amsthm,amssymb,amsmath} 
  \usepackage[alphabetic]{amsrefs}

  \usepackage[all]{xy}
  \usepackage{tikz-cd}
  \usepackage{quiver}

  \usepackage{bbm} 

  \usepackage[textsize=scriptsize]{todonotes}

  \usepackage{enumitem} 
  
    \usepackage[T1]{fontenc}
   \usepackage{kpfonts, baskervald}
   \usepackage{setspace}

  \theoremstyle{plain}
  \newtheorem{thm}{Theorem}[section]
  \newtheorem{prop}[thm]{Proposition}
  \newtheorem{cor}[thm]{Corollary}
  \newtheorem{lem}[thm]{Lemma}

  \theoremstyle{definition}
  \newtheorem{defi}[thm]{Definition}
   \newtheorem{rem}[thm]{Remark}
  \newtheorem{example}[thm]{Example}

\newcommand{\Cat}{\mathcal{C}at_\infty}
\renewcommand{\S}{\mathbb{S}}
\newcommand{\T}{\mathcal{S}}
\newcommand{\Th}{\mathrm{Th}}

\newcommand{\Fin}{\mathcal{F}in_\ast}
\newcommand{\Ass}{\mathcal{A}ss}
\newcommand{\LM}{\mathcal{LM}}
\newcommand{\E}{\mathbb{E}}
\newcommand{\Sp}{\mathcal{S}p}
\newcommand{\C}{\mathcal{C}}
\newcommand{\cC}{\mathcal{C}}

\newcommand{\dD}{\mathcal{D}}
\newcommand{\DD}{\mathcal{D}}

\newcommand{\Op}{\mathcal{O}}
\newcommand{\oO}{\mathcal{O}}
\newcommand{\OO}{\mathcal{O}}

\begin{document}

\onehalfspacing

\title{On Bialgebras, Comodules, Descent Data and Thom Spectra in $\infty$-categories}
\author{Jonathan Beardsley}
\maketitle

\begin{abstract}
This paper lays some of the foundations for working with not-necessarily-commutative bialgebras and their categories of comodules in $\infty$-categories. We prove that the categories of comodules and modules over a bialgebra always admit suitably structured monoidal structures in which the tensor product is taken in the ambient category (as opposed to a relative (co)tensor product over the underlying algebra or coalgebra of the bialgebra). We give two examples of higher coalgebraic structure: first, following Hess we show that for a map of $\E_n$-ring spectra $\phi\colon A\to B$, the associated $\infty$-category of descent data is equivalent to the category of comodules over $B\otimes_A B$, the so-called descent coring; secondly, we show that Thom spectra are canonically equipped with a highly structured comodule structure which is equivalent to the $\infty$-categorical Thom diagonal of Ando, Blumberg, Gepner, Hopkins and Rezk (which we describe explicitly) and that this highly structured diagonal decomposes the Thom isomorphism for an oriented Thom spectrum in the expected way.
\end{abstract}

\tableofcontents

\subsection*{Acknowledgements:}
It would not have been possible to write this paper without a great deal of help from many people, especially the following: Andrew Blumberg, Alexander Campbell, Tyler Lawson, Jack Morava, Denis Nardin, Eric Peterson, Maximilien P\'eroux, Maxime Ramzi, Emily Riehl, Liang Ze Wong, and Kirsten Zarek. I am especially grateful to Rune Haugseng, who has been immensely generous with his knowledge of the theory and practice of $\infty$-categories. Rune provided the author with sketch proofs of Theorem \ref{thm:monoidalEW} and Proposition \ref{prop:cocartlmodfibration} (obviously any mistakes or errors are the author's).  I am also in debt to the robust and flourishing online community of mathematicians working on homotopy theory and category theory. This work was partially supported by NSF grant DMS-1745583.  

\section{Introduction}

In this paper we develop some basic results about coalgebras and bialgebras in $\infty$-categories and give two applications in spectral algebra: one regarding $\infty$-categories of descent data, and another regarding coalgebra and comodule structures on Thom spectra. This last applications allows us to interpret (not necessarily commutative) Thom ring spectra from the point of view of non-commutative geometry.

In Section \ref{sec:coalgappendix} we make the relevant definitions (which have of course appeared in a number of other places) and ultimately show that if $H$ is a bialgebra with compatible $\E_k$-algebra and $\E_j$-coalgebra structures in a suitable $\infty$-category then its category of comodules has an $\E_k$-monoidal structure and its category of modules has an $\E_j$-monoidal structure (Theorem \ref{thm:bialgebramonoidalstructure}). We also make a detailed analysis of the relevant cocartesian morphisms of these monoidal structures to check that, given two left $H$-comodules $M$ and $N$, their tensor product is simply $M\otimes N$ (taken in the underlying category) and their comodule structure map is the same as the classical one, i.e.~the composite
\[
M\otimes N\to H\otimes M\otimes H\otimes N\overset{\sim}\to H\otimes H\otimes M\otimes N\to H\otimes M\otimes N
\]
of the respective coactions of $M$ and $N$ tensored together followed by a twist and the multiplication of $H$. Similarly, the tensor product of two left modules has structure map
\[
H\otimes M\otimes N\to H\otimes H\otimes M\otimes N\overset{\sim}\to H\otimes M\otimes H\otimes N\to M\otimes N
\]
given by the comultiplication of $H$ followed by the respective actions of $M$ and $N$.

Our first application, in Section \ref{sec:descendata}, is to give a description of the $\infty$-category of descent data for a map of $\E_n$-ring spectra $A\to B$ as the category of comodules over the so-called descent coring $B\otimes_AB$ (Theorem \ref{thm:descentdatacoring}). This seems to be a classical theorem, though the earliest reference the author has been able to find is \cite{brz}. Somewhat more recently it was proven in a model-category theoretic context by in \cite{hess}. In future work we will make a more in-depth study of $\infty$-categorical descent, and in particular how the above theorem manifests when $A\to B$ is a Galois or Hopf-Galois extension as defined in \cite{rog}. A result of possibly independent interest is a monoidal strengthening of the $\infty$-categorical Eilenberg-Watts theorem of \cite{ha}*{7.1.2.4} (Theorem \ref{thm:monoidalEW}).  

Our second application, in Sections \ref{sec:thomobjects} and \ref{sec:thomarecomods}, is to study the coalgebraic structure of Thom spectra. First, we show that any Thom spectrum $Th(f)$ over a ring spectrum $R$, constructed via a map $f\colon X\to Pic(R)$, is a comodule for the ``trivial'' Thom spectrum, $R\otimes \Sigma^\infty_+ X$, whose coalgebra structure is induced by the diagonal map of $X$ (Corollary \ref{thm:thomdiagonal}). While the existence of this ``Thom diagonal'' is a classical fact (cf. \cite{mahringthom} for one of many examples), our approach ensures that this coaction is fully homotopy coherent. In other words, it determines an object in an $\infty$-category of comodules, rather than just a comodule in the homotopy category. We also show that this diagonal is equivalent to the diagonal (implicitly) defined in \cite{abghr}, and that it plays a part in a Thom isomorphism whenever $\Th(f)$ is $E$-oriented for some $R$-algebra $E$ (Theorem \ref{thm:thomisomodules}). 

The role that our coherent diagonal plays in the Thom isomorphism is particularly interesting because, to the author's knowledge, no existing $\infty$-categorical reference gives a full decomposition of the Thom isomorphism into its classical components in which the first morphism is a \textit{structured} coaction, i.e.~the composite
\[E\otimes_R\Th(f)\to E\otimes_R \Th(f)\otimes  \Sigma^\infty_+X\to E\otimes_R E\otimes \Sigma^\infty_+X\to E\otimes \Sigma^\infty_+X\]
of the Thom diagonal, followed by the orientation, followed by the multiplication of $E$ (in the classical case this is described in \cite{lewis}). As a result, in the case that $f\colon X\to Pic(R)$ is $\E_k$-monoidal for some $k<n$ and use the fact that $\Th(f)$ is oriented with respect to itself, we are justified in saying that $\Th(f)$ is a non-commutative $\Sigma^\infty_+X$-torsor (or commutative, if everything in sight is $\E_\infty)$. Here we mean non-commutative in the sense of non-commutative algebraic geometry, in which one can mimic the classical notion of a $G$-torsor but replaces $G$ with a bialgebra (typically a Hopf algebra), schemes with rings, and pullbacks with tensor products. 

One shortcoming of the above analysis is that it does not fully expose the relationship between this torsor structure and the related descent theory. Classically, if $X\to S$ is a cover of schemes and a $G$-torsor for some group $G$, then one can perform ``descent along a torsor," (cf.~\cite{vistolistacks}*{Section 4.4}). In a homotopical context, for a ring map $\phi\colon A\to B$ with $G$ replaced by a spectral bialgebra $H$, this corresponds to a cosimplicial equivalence between the Amitsur complex for $\phi$ (i.e.~the Adams spectral sequence) and the $H$-cobar complex for $B$ over $A$. This equivalence is constructed in model categories of commutative and associative ring spectra in \cites{rog, froth} and is easy to construct by universal properties  when $\phi$ is a map of $\E_\infty$-ring spectra. We hope to give a purely $\infty$-categorical construction of this equivalence, and therefore a better understanding of non-commutative descent for Thom spectra, in future work. 

Finally in Appendix \ref{appendix: orientations}, we prove that (an $\infty$-categorical reformulation of) the classical notion of an $E$-orientation of a Thom spectrum, as described in \cite{lewis}, is indeed equivalent data the notion of orientation given in \cite{abghr}. This is well known to experts, and is implicit in \cite{acb}, but the author found it useful to explicitly describe the relationship between these two equivalent structures.

\section{Background}\label{sec:background}

Throughout, we use the $\infty$-categorical framework for higher category theory and derived algebra developed by Lurie in \cites{htt,ha}. For the most part our notation agrees with Lurie's. In particular, we use $\T$ for the $\infty$-category of spaces and $\Sp$ for the $\infty$-category of spectra. We will make heavy use of the theory of $\infty$-operads developed in \cite{ha}. In particular we will use:~the commutative $\infty$-operad $\Fin$; the little $n$-cubes operad $\E_n^\otimes$; the associative $\infty$-operad $\Ass^\otimes$, which is equivalent to $\E_1^\otimes$; and the $\infty$-operad $\LM^\otimes$ whose algebras are pairs $(A,M)$ comprising an algebra and a module over that algebra.

Given an $\E_n$-monoidal $\infty$-category we will write $Alg_{\E_k}(\C)$ for the category of $\E_k$-algebras in $\C$ for any $0\leq k\leq n$ and $LMod_A(\C)$ for the $\infty$-category of left modules over an $\E_1$-algebra in $\C$. Recall that if $LMod(\C)$ is the $\infty$-category of $\LM^\otimes$-algebras in $\C$ then $LMod_A(\C)$ can be obtained from $LMod(\C)$ as the following pullback of $\infty$-categories:
\[
\xymatrix{
LMod_A(\C)\ar[r]\ar[d] & LMod(\C)\ar[d] \\
\ast\ar[r]^{A} & Alg(\C)
}
\]
in which the right vertical map forgets the module and the bottom horizontal map picks out the algebra $A$.  More generally we will wish to consider categories of pairs $(A,M)$ where $A$ is an $\E_k$-algebra and $M$ is a left module over the underlying $\E_1$-algebra of $A$. 

\begin{defi}\label{def:Ekmodules}
Let $LMod_{\E_k}(\C)$ be the $\infty$-category defined by the following pullback of $\infty$-categories:
\[
\xymatrix{
LMod_{\E_k}(\C)\ar[r]\ar[d] & LMod(\C)\ar[d]\\
Alg_{\E_k}(\C)\ar[r] & Alg(\C)
}
\]
in which the bottom horizontal map takes an $\E_k$-algebra to its underlying $\E_1$-algebra. The objects of $LMod_{\E_k}(\C)$ can be thought of as pairs $(A,M)$ where $A$ is an $\E_k$-algebra of $\C$ and $M$ is a left module over the underlying $\E_1$-algebra of $A$. 
\end{defi}

\begin{rem}\label{rem:Ekforget}
Recall that objects of $LMod(\C)$ are commutative triangles
\[
\xymatrix{
\LM^\otimes\ar[r]^F\ar[dr] & \C^\otimes\times_{\E_n^\otimes}\E_1^\otimes\ar[d]\\
& \E_1^{\otimes}
}
\]
in which the horizontal functor $F$ takes inert morphisms of $\LM^\otimes$ to inert morphisms of $\C^\otimes\times_{\E_n^\otimes}\E_1^\otimes$ (with respect to the cocartesian fibration $\C^{\otimes}\times_{\E_n^\otimes}\E_1^\otimes\to\E_1^\otimes$). Recall also that the objects of $\LM^\otimes_{\langle 1\rangle}$ are $\mathfrak{a}$ and $\mathfrak{m}$, corresponding to an algebra and a module respectively. By precomposing with the inclusion $\mathfrak{m}\colon \ast\to \LM^\otimes$ we obtain a forgetful functor $LMod(\C)\to \C$ which takes a pair $(A,M)$ to $M$ as an object of $\C$. By further precomposing with the projection $LMod_{\E_k}(\C)\to LMod(\C)$ we obtain a similar forgetful functor $LMod_{\E_k}(\C)\to \C$.
\end{rem}

\begin{rem}
 Note that, by pasting pullback diagrams as follows,
	
	\[
	\xymatrix{
		LMod_A(\C)\ar[r]\ar[d] & LMod_{\E_k}(\C)\ar[r]\ar[d] & LMod(\C)\ar[d]\\
		\ast\ar[r]^A & Alg_{\E_k}(\C)\ar[r] & Alg(\C)
	}
	\]
	we have that the pullback of the cartesian fibration $LMod_{\E_k}(\C)\to Alg_{\E_k}(\C)$ along the morphism $\ast\to Alg_{\E_k}(\C)$ is equivalent to $LMod_A(\C)$. As a result, we will not differentiate between the two $\infty$-categories. The same goes for $LCoMod_A(\C)$ for an $\E_k$-coalgebra $A$.
\end{rem}

\begin{prop}\label{prop:LModEkCartesianfibration}
The projection $LMod_{\E_k}(\C)\to Alg_{\E_k}(\C)$ of Definition \ref{def:Ekmodules} is a cartesian fibration of $\infty$-categories with respect to which a morphism $(A,M)\to (B,N)$ is cartesian if and only if it induces an equivalence after applying the forgetful functor $LMod_{\E_k}(\C)\to\C$ of Remark \ref{rem:Ekforget}. 
\end{prop}

\begin{proof}
From \cite{ha}*{4.2.3.2}, we have that $LMod(\C)\to Alg(\C)$ is a cartesian fibration and that a morphism of $LMod(\C)$ is cartesian if and only if it is an equivalence after application of the forgetful functor $LMod(\C)\to \C$ described in Remark \ref{rem:Ekforget}. The fact that $LMod_{\E_k}(\C)\to Alg_{\E_k}(\C)$ is cartesian then follows from the fact that cartesian fibrations are stable under pullback (cf.~\cite{htt}*{2.4.2.3}). Additionally, by \cite{htt}*{2.4.1.12, 2.4.2.8}, a morphism in $LMod_{\E_k}(\C)$ is cartesian over $Alg_{\E_k}(\C)$ if and only if its image under the projection $LMod_{\E_k}(\C)\to LMod(\C)$ is also cartesian. As a result, the characterization of cartesian edges in $LMod(\C)$ extends to a characterization of cartesian edges in $LMod_{\E_k}(\C)$. 
\end{proof}

\begin{rem}
One can deduce, for instance by applying \cite{kerodon}*{5.5.4.13} and the straightening-unstraightening equivalence of \cite{htt}, that the cartesian morphism with target $(B,N)$ in either $LMod(\C)$ or $LMod_{\E_k}(\C)$ and over $\phi\colon A\to B$ in either $Alg(\C)$ or $Alg_{\E_k}(\C)$ can be thought of as $(A,\phi^\ast(N))\to (B, N)$ given by the pair $(\phi, id_N)$ where the $A$-action on $\phi^\ast(N)\simeq N$ is given by $\phi$ followed by the $B$-module structure on $N$, i.e.~$A\otimes N\to B\otimes N\to N$. 
\end{rem}

\begin{defi}
	Let $\Op^\otimes$ be an $\infty$-operad and $\C$ an $\Op$-monoidal $\infty$-category defined by a cocartesian fibration $p\colon\C^{\otimes}\to\Op^\otimes$ as in \cite{ha}*{2.1.1.10}. Then the \textit{dual fibration} of $p$ is the cartesian fibration $p^\vee\colon (\C^\otimes)^\vee\to\Op^\otimes$ as defined in \cite{bgn}. We define the \textit{opposite $\Op$-monoidal structure on $\C$} to be the $\Op$-monoidal structure on $\C$ determined by the opposite of the dual fibration, i.e.~the cocartesian fibration $(p^\vee)^{op}\colon ((\C^\otimes)^\vee)^{op}\to \Op^\otimes$. We will call the fibration $(p^\vee)^{op}$ the ``fiberwise opposite'' of $p$. To simplify notation we will follow \cite{haugsengmates} and write $p_{op}\colon \C_\otimes^{op}\to \Op^\otimes$ for the fiberwise opposite of $p$.
\end{defi}

\begin{rem}
	It follows from \cite{bgn}*{Theorem 1.7} that if $p\colon\C^{\otimes}\to\Op^\otimes$ is an $\Op$-monoidal structure on $\C$ then the functor $\Op^\otimes\to\Cat$ corresponding to the fibration $p_{op}\colon \C_{\otimes}^{op}\to \Op^\otimes$ is equivalent to the composition $op\circ St(p)\colon \Op^\otimes\to\Cat$, where $St(-)$ is the straightening functor of \cite{htt}. In other words, for any $X\in\Op^\otimes$, the fiber $(\C_\otimes^{op})_X$ is equivalent to $(\C^\otimes_X)^{op}$. 
\end{rem}

\begin{defi}\label{defi:oppositefunctor}
Let $\C$ and $\dD$ be $\Op$-monoidal $\infty$-categories and $f\colon \C\to \dD$ be an $\Op$-monoidal functor. Then by functoriality of taking the dual fibration and taking opposites there is an induced $\oO$-monoidal functor which we denote by $f_{op}\colon \C^{op}\to \dD^{op}$. 
\end{defi}

\begin{rem}
Going forward, when there is no chance of ambiguity, we will not reference the defining cocartesian fibration $p\colon \C^\otimes\to \Op^\otimes$ for an $\Op$-monoidal $\infty$-category, nor that of the opposite monoidal structure. Whenever we fix an $\Op$-monoidal $\infty$-category $\C$ we will implicitly be fixing a cocartesian fibration $p\colon \C^\otimes\to \Op^\otimes$ and implicitly endowing the opposite category $\C^{op}$ with the $\Op$-monoidal structure determined by $p_{op}$. 
\end{rem}

\section{Coalgebra}\label{sec:coalgappendix}

\begin{defi}\label{coalgs}
	Let $\C$ be an $\mathcal{O}$-monoidal $\infty$-category for $\mathcal{O}^\otimes$ an $\infty$-operad. Then define the $\infty$-category of $\Op$-coalgebras in $\C$ to be $(Alg_{\Op}(\C^{op}))^{op}$, which we will denote by $CoAlg_{\Op}(\C)$. If ${\Op}^\otimes\simeq\E_1^\otimes$ we will write $CoAlg(\C)$ for the $\infty$-category of coassociative coalgebras in $\C$. 
\end{defi}

\begin{defi}[Comodules]\label{comods}
	Let $\C$ be an $\E_n$-monoidal $\infty$-category and let $A$ be an object of $Alg_{\E_k}(\C^{op})$ for $0<k\leq n$. Then there is an $\infty$-category of left $A$-modules $LMod_A(\C^{op})$.  We define the category of left comodules over $A$ to be the $\infty$-category $LMod_A(C^{op})^{op}$. We will denote this category by $LCoMod_A(\C)$. Similarly  for any $k\leq n$ we define $LCoMod_{\E_k}(\C)$ to be $LMod_{\E_k}(\C^{op})^{op}$.
\end{defi}

\begin{example}\label{prop:spacecoalgs}
	By \cite{beardsleyperoux}*{2.2}, if $\C$ is a cartesian symmetric monoidal $\infty$-category, i.e.~one in which the tensor product of two objects is given by their cartesian product, then every object is a canonically a cocommutative coalgebra. Moreover for a fixed object $X\in \C$ there is an equivalence of $\infty$-categories between $LCoMod_X(\C)$ and the slice category $\C_{/X}$. In particular, every space is a cocommutative coalgebra of $\T$ and a map of spaces $X\to Y$ induces a $Y$-comodule structure on $X$ via the composite $X\to X\times X\to X\times Y$.
\end{example}

\begin{rem}\label{expo:comodulestructuremap}
Definition \ref{comods} is somewhat abstract given that the monoidal structure is on $\C^{op}$ is determined by the opposite of the dual of the fibration determining the monoidal structure on $\C$, so we will describe part of it more explicitly in the case of $\OO=Ass^{\otimes}$. Let $p\colon \C^{\otimes}\to Ass^{\otimes}$ be the fibration determining the monoidal structure on $\C$ and $p_{op}\colon C_\otimes^{op}\to Ass^\otimes$ the fiberwise opposite. By \cite{bgn}*{1.6, 3.5} we know that a morphism $(\phi,\alpha)\colon x\to y$ in $C_\otimes^{op}$ is a cospan $x\overset{\phi}\to u\overset{\alpha}\leftarrow y$ in $\C^{\otimes}$ where $p(\phi)\simeq p_{op}(\phi,\alpha)$, $\phi$ is cocartesian in $\C^\otimes$, and $p(\alpha)$ is equivalent to the identity. Moreover, the morphism $(\phi,\alpha)$ is cocartesian in $\C_\otimes^{op}$ exactly when $\alpha$ is an equivalence in $\C^{\otimes}$. 

Now suppose we are given a morphism of $\infty$-operads over $Ass^\otimes$:
\[
\xymatrix{
\LM^\otimes\ar[r]^F\ar[dr] & (\C_\otimes)^{op}\ar[d]\\
& Ass^\otimes
}
\]
determining an algebra and a module over it in $\C^{op}$. If we say that $F(\mathfrak{a})=A$ and $F(\mathfrak{m})=M$ then the morphism of $\infty$-operads $F$ determines, among other things, a commutative triangle in $\C_{\otimes}^{op}$ of the following form:
\[
\xymatrix{
F((\mathfrak{a},\mathfrak{m}))\simeq (A,M)\ar[d]_-{\delta}\ar[r]^-{\beta} & A\otimes M\ar[dl]^-{\gamma}\\
F(\mathfrak{m})\simeq M & 
}
\]
where $\beta$ is cocartesian, $\delta$ is the image of the active map $(\mathfrak{a},\mathfrak{m})\to \mathfrak{m}$ in $\LM^\otimes$ under $F$, and $\gamma$ is induced by the fact that $\beta$ is cocartesian. Here, $\gamma$ is the map witnessing the action of $A$ on $M$ in $\C^{op}$ (ignoring higher coherences). By replacing each of $\beta$, $\gamma$ and $\delta$ with their corresponding cospans in $\C^\otimes$, we have a diagram:
\[
\xymatrix{
(A,M)\ar[r]^{\beta_L} \ar[d]_{\delta_L} & U & A\otimes M\ar[l]_{\beta_R}\ar[dl]^{\gamma_L} \\
V & W & \\
M\ar[u]^{\delta_R}\ar[ur]_{\gamma_R} & &
}
\]
in which $\beta_L$, $\delta_L$ and $\gamma_L$ are cocartesian in $\C^\otimes$. Moreover, $\beta_L$ and $\delta_L$ project to the active map $\langle 2\rangle\to\langle 1\rangle$, $\gamma_L$ projects to the identity map $\langle 1\rangle\to \langle 1\rangle$ and  $\delta_R$ and $\gamma_R$ project to degenerate morphisms in $Ass^\otimes$ under $p$. Finally, because $\beta$ is cocartesian, $\beta_R$ is an equivalence. Since $\beta_L$ and $\delta_L$ are both cocartesian in $\C^\otimes$ they must both be equivalent to the tensor product structure map $(A,M)\to A\otimes M$ in $\C^\otimes$. Similarly by cocartesianness, $\gamma_L$ must be an equivalence $A\otimes M\to A\otimes M$. Now note that this diagram must commute in $\C_\otimes^{op}$ and that composition in $\C_\otimes^{op}$ is given by taking pushouts. Therefore there is a commutative diagram in $\C^\otimes$:

\[\begin{tikzcd}
	{(A,M)} && {A\otimes M} && M \\
	& {A\otimes M} && {A\otimes M} \\
	&& {A\otimes M}
	\arrow["{\beta_L}", from=1-1, to=2-2]
	\arrow[from=2-2, to=3-3]
	\arrow["{\beta_R}"', from=1-3, to=2-2]
	\arrow["{\gamma_L}", from=1-3, to=2-4]
	\arrow[from=2-4, to=3-3]
	\arrow["{\gamma_R}"', from=1-5, to=2-4]
	\arrow[from=1-3, to=2-2]
	\arrow[from=1-3, to=2-4]
	\arrow["{\delta_L}"', curve={height=30pt}, from=1-1, to=3-3]
	\arrow["{\delta_R}", curve={height=-30pt}, from=1-5, to=3-3]
	\arrow["\lrcorner"{anchor=center, pos=0.125, rotate=135}, draw=none, from=3-3, to=1-3]
\end{tikzcd}\]
in which the central square is a pushout square each of whose sides is an equivalence. It follows that $\delta_R\simeq \gamma_R$ and $\delta_L\simeq \beta_L$. In particular, the morphism $\gamma_R$ is precisely the ``coaction map'' of $A$ on $M$ in $\C$. This is the way in which the specific morphism in $\C$ describing the coaction may be extracted from the $\infty$-operadic data.

\end{rem}

\subsection{Bialgebras}

By applying \cite{ha}*{3.2.4.3} to the bifunctor $\E_j^\otimes\times\E^\otimes_k\to\E^\otimes_{k+j}$ of \cite{ha}*{5.1.2.1} we have that if the tensor product in an $\E_{k+j}$-monoidal $\infty$-category preserves geometric realizations then the $\infty$-category of $\E_k$-algebras in an $\C$ is always at least $\E_j$-monoidal. However, in general, $Alg_{\E_k}(\C)$ will not be $\E_\ell$-monoidal for $\ell>j$. As a result, if we are interested in discussing bialgebras in an $\E_n$-monoidal $\infty$-category, our constructions only allow us to work with bialgebras that have an $\E_j$-comonoidal structure and an $\E_k$-monoidal structure for $j,k\geq 0$ and $j+k\leq n$. We will call such bialgebras $(j,k)$-bialgebras.  

\begin{defi}[$(j,k)$-Bialgebras]\label{cohopfalgs}
	Let $\C$ be an $\E_n$-monoidal $\infty$-category. Then for any $j\leq n$ there is an $\E_{n-j}$-monoidal $\infty$-category of $\E_j$-coalgebras in $\C$, $CoAlg_{\E_j}(\C)$. For each $k\leq n-j$, there is an $\infty$-category $Alg_{\E_k}(CoAlg_{\E_j}(\C))$. For a fixed $j,k<n$, we call $Alg_{\E_k}(CoAlg_{\E_j}(\C))$ the category of $(j,k)$-bialgebras in $\C$. We will denote this category by $BiAlg_{j,k}(\C)$ where the first index gives the ``degree'' of cocommutativity, and the second index gives the ``degree'' of commutativity. 
\end{defi}

\begin{rem}
	In the above definition, if $n=\infty$, so that $\C$ and $\C^{op}$ are symmetric monoidal, then $Alg_{\E_k}(\C^{op})$ is again symmetric monoidal (cf.~\cite{ha}*{3.2.4.4}). As such, in a symmetric monoidal $\infty$-category, we can define $BiAlg_{j,k}(\C)$ for arbitrary $j$ and $k$. 
\end{rem}

\begin{rem}
	Note that for an $\E_{k+j}$-monoidal $\infty$-category $\C$, an object $A$ of $Alg_{\E_j}(\C^{op})$ is a section of the cocartesian fibration $\C_{op}^\otimes\to \E_j^\otimes$ defining the $\E_j$-monoidal structure on $\C^{op}$. Therefore the image of the inclusion of the base point $\{\ast\}\to\langle 1\rangle$ in $\E_j^\otimes$ (cf.~\cite{ha}*{5.1.0.2}) induces an algebra unit map $1_\C\to A$ in $\C^{op}$. Hence $A$ is equipped with a counit $\varepsilon\colon A\to 1_\C$ in $\C$. Similarly, $A$ is equipped with a comultiplication $\Delta\colon A\to A\otimes A$ which is homotopy cocommutative when $j>1$. This map corresponds to (the opposite of) a cocartesian lift of a rectilinear mapping of a pair of $j$-cubes into a $j$-cube. 
\end{rem}

\begin{rem}
	Recall that when defining an affine commutative monoid scheme, one defines it to be a commutative monoid object in the category of affine schemes. As a result, an affine commutative monoid scheme is both a commutative monoid and a cocommutative comonoid (via the diagonal map of the underlying scheme, or $Spec(-)$ of the multiplication), and more importantly, these two structures are compatible. In other words, to produce a bialgebra, we can either specify an algebra whose structure maps are maps of coalgebras, or a coalgebra whose structure maps are maps of algebras. Both of these conditions will produce the necessary compatibility between these structures. In a 1-category this structure can be encoded by certain diagrams or, equivalently, the symmetric monoidal theory (or PROP) of bialgebras. The machinery of \cite{ha} and \cite{htt} does not immediately lend itself to such a description. There does exist a theory of $\infty$-properads which can be used to parameterize compatible homotopy coherent algebra and coalgebra structures described in \cite{hry}, but it is very difficult to make that machinery compatible with the rest of the framework for derived algebra described in \cite{htt,ha}. Luckily, however, we will not need such a compact description of bialgebras in this paper.
\end{rem}

\begin{prop}\label{prop:functorofbialgebra}
	Let $\C$ and $\dD$ be $\E_n$-monoidal $\infty$-categories and $f\colon\C\to\dD$ an $\E_n$-monoidal functor. Then for $j+k=n$, if $H$ is a $(j,k)$-bialgebra in $\C$ then $f(H)$ is a $(j,k)$-bialgebra in $\dD$. 
\end{prop}

\begin{proof}

Being $\E_n$-monoidal, $f$ induces an $\E_n$-monoidal functor $f_{op}\colon \cC^{op}\to\DD^{op}$. As such, $f_{op}$ induces an $\E_k$-monoidal functor $Alg(f_{op})\colon Alg_{\E_j}(\C^{op})\to Alg_{\E_j}(\dD^{op})$. This in turn induces an $\E_k$-monoidal functor $CoAlg(f)\colon CoAlg_{\E_j}(\C)\to CoAlg_{\E_j}(\dD)$ which, being $\E_k$-monoidal, preserves $\E_k$-algebras in $CoAlg_{\E_j}(\C)$ and therefore preserves $(j,k)$-bialgebras. 

\end{proof}

\subsection{Tensor Products of Modules and Comodules over Bialgebras}

In general, categories of left (or right) comodules over an $\E_n$-coalgebra in an $\E_n$-monoidal $\infty$-category $\C$ do not necessarily admit any monoidal structure (not even for $n=\infty$). This is true even in the discrete case, so we do not expect it to be true in the $\infty$-categorical case. Approaching the problem na\"ively one might expect that, for a coalgebra $A\in CoAlg_{\E_n}(\C)$, i.e.~$A\in Alg_{\E_n}(\C^{op})$, we could follow Lurie in \cite{ha}*{5.1.4} and apply \cite{ha}*{4.8.5.20} to obtain an $\E_{n-1}$-monoidal structure on $LMod_A(\C^{op})$. Unfortunately \cite{ha}*{4.8.5.20} in not applicable in this case unless $\C$ has the unlikely property that its tensor product commutes with totalizations of cosimplicial objects (as the application of ibid.~requires the tensor product of $\C^{op}$ to commute with geometric realizations of simplicial objects). 

 Luckily, if $A$ is not just a coalgebra but a bialgebra, we can endow $LCoMod_A(\C)$ with a monoidal structure (cf.~Theorem \ref{thm:bialgebramonoidalstructure}). Note that this monoidal structure is \textit{not} the ``relative tensor product'' of \cite{ha}*{4.4}, which is, by the preceding paragraph, often impossible to construct. Instead, the monoidal structure on $LCoMod_A(\C)$ is, after forgetting to $\C$, equivalent to the ambient tensor product in $\C$ itself. So for instance the unit of this monoidal structure on $LCoMod_A(\C)$ is $1_\C$ rather than $A$. A similar statement applies to the case of left modules over a bialgebra.

We begin with some preliminaries:

\begin{defi}
	A bifunctor of $\infty$-operads $\mathcal{P}^\otimes\times \mathcal{Q}^\otimes\to \mathcal{O}^\otimes$, as defined in \cite{ha}*{2.2.5.3},  is a functor of $\infty$-categories that is compatible with the smash product functor of finite pointed sets $Fin_\ast\times Fin_\ast\overset{\wedge}\to Fin_\ast$ via the defining fibrations of the $\infty$-operads in question. In other words, there is a commutative diagram of $\infty$-categories
	\[
	\xymatrix{
		\mathcal{P}^\otimes\times\mathcal{Q}^\otimes\ar[d]\ar[r] & \mathcal{O}^\otimes\ar[d]\\
		Fin_\ast\times Fin_\ast\ar[r]^-{\wedge} & Fin_\ast
	}
	\] 
	where the vertical maps are the ones defining the $\infty$-operad structures of $\mathcal{P}$, $\mathcal{Q}$ and $\mathcal{O}$.
\end{defi}

\begin{rem}\label{rem:additivitybifunctor}
	Recall that there is an important bifunctor $\E_n^\otimes \times \E_k^\otimes \to \E_{n+k}^\otimes $ given, essentially, by taking products of cube embeddings, as described in \cite{ha}*{5.1.2.1}. This is the bifunctor which ultimately leads to the $\infty$-categorical version of Dunn additivity (cf.~\cite{ha}*{5.1.2.2}). We will sometimes use the obvious extension of this bifunctor to one of the form $\E_n^\otimes\times\E_j^\otimes\to \E_{n+k}$ for any $0\leq j\leq k$
\end{rem}

\begin{lem}\label{lem:LMbifunc}
	For any $k\leq n-1$ there is a bifunctor of $\infty$-operads $\LM^\otimes\times\E_{k}^\otimes\to \E_{n}^\otimes$.
\end{lem}

\begin{proof}
	The lemma follows from considering the following commutative diagram:
	
	\[
	\xymatrix{
		\LM^\otimes\times\E_k^\otimes\ar[r]\ar[d] &  \E_1^\otimes\times\E_k^\otimes\ar[d]\ar[r] & \E_{k+1}\ar[d]\ar[r] & \E_n^\otimes\ar[d] \\
		Fin_\ast\times Fin_\ast\ar[r]^-{id} &  Fin_\ast\times Fin_\ast\ar[r]^-{\wedge} & Fin_\ast\ar[r]^-{id} & Fin_\ast
	}
	\]
	
	The first and third squares commute because there are functors of $\infty$-operads $\LM^\otimes\to \E_1^\otimes$ and $\E_{k+1}^\otimes\to \E_n^\otimes$ by \cite{ha}*{4.2.1.9, 5.1.1.5}. The second square commutes because it is the bifunctor of Remark \ref{rem:additivitybifunctor}. 
\end{proof}

\begin{rem}\label{expose:leftmodulemonoidal}
	Throughout this section, especially in Theorem \ref{thm:bialgebramonoidalstructure}, we will make heavy use of \cite{ha}*{3.2.4.3} which states that, given an $\OO$-monoidal $\infty$-category $p\colon \C^\otimes\to\OO^\otimes$ and a bifunctor of $\infty$-operads $\mathcal{P}^\otimes\times \mathcal{Q}^\otimes\to \OO^\otimes$, there is a $\mathcal{Q}$-monoidal structure on the $\infty$-category of $\mathcal{P}$-algebras in $\C$ (with respect to the $\OO$-monoidal structure), i.e.~a cocartesian fibration $q\colon Alg_{\mathcal{P}/\OO}(\C)^\otimes\to\mathcal{Q}^\otimes$. 
	
	 Taking $\mathcal{P}^\otimes=\LM^\otimes$, $\mathcal{Q}^\otimes=\E_k^\otimes$ and $\OO^\otimes=\E_n^\otimes$, for $k\leq n-1$ as in Lemma \ref{lem:LMbifunc}, we have an $\E_k$-monoidal structure on $Alg_{\LM/\E_n}(\C)\simeq Alg_{\LM/\E_1}(\C)\simeq LMod(\C)$, where the objects of the right hand side can be thought of as pairs $(A,M)$, with $A\in Alg(\C)$ and $M$ a left $A$-module. Going forward we will always use the notation $LMod(\C)$ instead of $Alg_{\LM/\E_n}(\C)$. 
By using the bifunctor $\E_k^\otimes\times \E_1^\otimes\to \E_n^\otimes$ in the same way (for $k+j\leq n$), we obtain an $\E_k$-monoidal structure on $Alg_{\E_1}(\C)$.
	
	From \cite{ha}*{3.2.4.3 (4)}, we have that a morphism $\alpha\in LMod(\C)^\otimes$ (resp. $\alpha\in Alg_{\E_j}(\C)^\otimes$), i.e.~a certain kind of natural transformation $F\to G$ between functors $\LM^\otimes\to\C^{\otimes}$ (resp. functors $\E_j^\otimes\to\C^\otimes$), is cocartesian if and only if for each object $X\in \LM^\otimes$ (resp. $X\in \E_j^\otimes$) the morphism $\alpha(X)\colon F(X)\to G(X)$ is cocartesian in $\C^\otimes$. As a result, the ``forgetful'' functors $LMod(\C)^\otimes\to Alg_{\E_1}(\C)^\otimes$ and $LMod(\C)^\otimes\to \C^\otimes\times_{\E^\otimes_n}\E^\otimes_k$, given by precomposition with $Ass^\otimes\hookrightarrow \LM^\otimes$ and evaluation at $\mathfrak{m}$ respectively, both preserve cocartesian morphisms over $\E_k^\otimes$ and are therefore $\E_k$-monoidal. 
	
\end{rem}

\begin{lem}\label{lemma:monoidalstructureonLMod}
If $\C$ is an $\E_n$-monoidal $\infty$-category for $n>1$ then $LMod_{\E_k}(\C)$ and $Alg_{\E_k}(\C)$ both admit $\E_{n-k}$-monoidal structures for $0<k<n$. Moreover for each $0<k<n$ the projections $LMod_{\E_k}(\C)\to Alg_{\E_k}(\C)$ and $LMod_{\E_k}(\C)\to \C$ are $\E_{n-k}$-monoidal.
\end{lem}

\begin{proof}

Fix some $0<k<n$. Then $LMod(\C)$, $Alg(\C)$ and $Alg_{\E_k}(\C)$ are all $\E_{n-k}$-monoidal categories.  For the first two categories this follows from Remark \ref{expose:leftmodulemonoidal} and from using the bifunctor of \ref{rem:additivitybifunctor} for the third. Moreover, because all of these monoidal structures have cocartesian morphisms which are determined by the cocartesian morphisms of $\C^\otimes$ (again by \cite{ha}*{3.2.4.3 (4)}), the forgetul functors $LMod(\C)\to Alg(\C)$ and $Alg_{\E_k}(\C)\to Alg(\C)$ are both $\E_{n-k}$-monoidal over $\C$. Therefore we can take the pullback of the cospan $Alg_{\E_k}(\C)\to Alg(\C)\leftarrow LMod(\C)$ in the slice category of $\E_{n-k}$-monoidal $\infty$-categories over $\C$. 
This pullback presents $LMod_{\E_k}(\C)$ as an $\E_{n-k}$-monoidal $\infty$-category such that the projections $LMod_{\E_k}(\C)\to \C$ and  $LMod_{\E_k}(\C)\to Alg_{\E_k}(\C)$ are $\E_{n-k}$-monoidal. 
\end{proof}

\begin{rem}
Note that by taking fiberwise opposites appropriately, Lemma \ref{lemma:monoidalstructureonLMod} implies that $CoAlg_{\E_k}(\C)$ and $LCoMod_{\E_k}(\C)$ also admit $\E_{n-k}$-monoidal structures.
\end{rem}

\begin{prop}\label{prop:cocartlmodfibration}
Let $\C$ be an $\E_n$-monoidal $\infty$-category. Then there is a commutative triangle 
\[
\xymatrix{
LCoMod_{\E_k}(\C)^\otimes\ar[r]^f\ar[dr]_p & CoAlg_{\E_k}(\C)^\otimes\ar[d]^q\\
& \E_{n-k}^\otimes
}
\] where $p$ and $q$ and $f$ are all cocartesian fibrations and $f$ is the opposite of the $\E_{n-k}$-monoidal projection functor $LMod_{\E_k}(\C^{op})^\otimes\to Alg_{\E_k}(\C^{op})^\otimes$ constructed in Lemma \ref{lemma:monoidalstructureonLMod}.
\end{prop}

\begin{proof}
By Lemma \ref{lemma:monoidalstructureonLMod} and Proposition \ref{prop:LModEkCartesianfibration} there is an $\E_{n-k}$-monoidal projection functor $F\colon LMod_{\E_k}(\C^{op})^\otimes\to Alg_{\E_k}(\C^{op})^\otimes$ which is a fiberwise cartesian fibration. By taking fiberwise opposites we obtain a commutative triangle 
\[
\xymatrix{
LCoMod_{\E_k}(\C)^\otimes\ar[r]^{F_{op}}\ar[dr]_p & CoAlg_{\E_k}(\C)^\otimes\ar[d]^q\\
& \E_{n-k}^\otimes
}
\]
in which $p$ and $q$ are cocartesian and $F_{op}$ is cocartesian when restricted to fibers over $\langle m\rangle\in\E_{n-k}^\otimes$. By taking $f= F_{op}$ it only remains to show that $f$ is in fact a cocartesian fibration rather than only cocartesian on fibers. We use \cite{haugseng_melani_safronov}*{A.1.8}, which gives sufficient conditions for a fiberwise cocartesian fibration to be cocartesian. We will check the three conditions of ibid.~to complete the proof. Conditions (1) and (3) are immediately satisfied by construction. Condition (2) is satisfied because $F_{op}$ is $\E_{n-k}$-monoidal. For condition (4) we need to check that the for each morphism $\phi\colon \langle \ell\rangle\to\langle m\rangle$ in $\E_{n-k}^\otimes$ the induced functor on fibers $\phi_!\colon LCoMod_{\E_k}(\C)^{\ell}\to LCoMod_{\E_k}(\C)^m$ takes cocartesian lifts of morphisms of $CoAlg_{\E_k}(\C)^{\ell}$ to cocartesian lifts of morphisms of $CoAlg_{\E_k}(\C)^m$ (i.e.~if $\alpha$ is a $f_\ell$-cocartesian then $\phi_!(\alpha)$ is $f_m$-cocartesian).  By Proposition \ref{prop:LModEkCartesianfibration}, a cocartesian lift of a morphism in $CoAlg_{\E_k}(\C)^\ell$ to $LCoMod_{\E_k}(\C)^\ell$ is one which is an equivalence when projected down to $\C^\ell$. It is straightforward to check that $\phi_!$ will preserve this structure whenever $\phi$ is either inert or active in $\E_{n-k}^\otimes$. By \cite{ha}*{2.1.2.4} we have that every morphism in an $\infty$-operad factors as an inert morphism followed by an active morphism. Therefore, by factoring $\phi$ in such a way, the proposition is proven.
\end{proof} 

\begin{rem}\label{expose:cocartmorphsincomod}
\sloppy
It will be useful going forward to have explicit descriptions of the morphisms in $LCoMod_{\E_k}(\C)$ which are cocartesian with respect to the functor $f$ defined in Proposition \ref{prop:cocartlmodfibration}. For any morphism in $CoAlg_{\E_k}(\C)^\otimes$ contained in a single fiber over $\langle \ell\rangle\in\E_{n-k}$ the cocartesian lift to $LCoMod_{\E_k}(\C)^\otimes$ lies in $LCoMod_{\E_k}(\C)^\ell$ and is just the juxtaposition of the $\ell$ cocartesian lifts obtained over each coordinate. In general however we can construct cocartesian lifts of any morphism of $CoAlg_{\E_k}(\C)$ by tracing through the proof of \cite{haugseng_melani_safronov}*{A.1.8}. In  what follows we will use the same notation of $p$, $q$ and $f$ from the statement of Proposition \ref{prop:cocartlmodfibration}. 

Let $\hat{M}$ be an object of $LCoMod_{\E_k}(\C)^\otimes$ with $f(\hat{M})\simeq\hat{A}\in CoAlg_{\E_k}(\C)^\otimes$ ($\hat{M}$ can be thought of as a finite list of pairs $(A,M)$ where $A$ is a coalgebra and $M$ is an $A$-comodule, and $\hat{A}$ is the associated list of coalgebras). Let $\delta\colon \hat{A}\to \hat{B}$ be a morphism in $CoAlg_{\E_k}(\C)^{\otimes}$ with $q(\delta)\simeq(\gamma\colon \langle \ell\rangle\to\langle m\rangle)$ in $\E_{n-k}^\otimes$. First notice that because $p$ is a cocartesian fibration there is a cocartesian morphism $\beta\colon \hat{M}\to\gamma_!\hat{M}$ in $LCoMod_{\E_k}(\C)^\otimes$ lifting $\gamma$. Because $f$ is $\E_{n-k}$-monoidal (i.e.~it takes $p$-cocartesian morphisms to $q$-cocartesian morphisms) we have a $q$-cocartesian morphism $f(\beta)\colon f(\hat{M})\simeq\hat{A}\to f(\gamma_!\hat{M})$ again lifting $\gamma$. Because $f(\beta)$ is $q$-cocartesian there is a factorization
\[
\xymatrix{
\hat{A}\ar[r]^{f(\beta)}\ar[dr]_{\delta} & f(\gamma_!(\hat{M}))\ar[d]^\alpha\\
& \hat{B}
}
\] 
in which $\alpha$ is a morphism in the fiber of $CoAlg_{\E_k}(\C)^\otimes$ over $\langle m\rangle$. Because $\alpha$ is a morphism in a single fiber and $f$ is a fiberwise cocartesian fibration, we can lift it to a cocartesian morphism $\epsilon\colon \gamma_!\hat{M}\to\alpha_!\gamma_!\hat{M}$. The above construction follows the proof of  \cite{haugseng_melani_safronov}*{A.1.8} verbatim (with some variable and category names changed) and so it follows that the composition \[\epsilon\circ \beta\colon \hat{M}\to\gamma_!\hat{M}\to\alpha_!\gamma_!\hat{M}\] is a cocartesian lift of $\delta$ to $LCoMod_{\E_k}(\C)$. 

Now we use the above formula to determine cocartesian lifts of certain morphisms in $CoAlg_{\E_k}(\C)$ which will be of interest to us. Suppose that $A_1\to B_1$ and $A_2\to B_2$ are morphisms in $CoAlg_{\E_k}(\C)$. Then there is a morphism $(A_1,A_2)\to B_1\otimes B_2$, namely the composition $(A_1,A_2)\to (B_1,B_2)\to B_1\otimes B_2$. The cocartesian morphism $(A_1,A_2)\to A_1\otimes A_2$ then induces a unique morphism $\psi\colon A_1\otimes A_2\to B_1\otimes B_2$ (the ``obvious'' tensor product of the maps $A_1\to B_1$ and $A_2\to B_2$).  The description of $f$-cocartesian lifts in the preceding paragraphs, along with the fact that the projection maps are $\E_{n-k}$-monoidal as shown in Lemma \ref{lemma:monoidalstructureonLMod},  implies that the cocartesian lift of the morphism $(A_1,A_2)\to B_1\otimes B_2$ with domain $((A_1,M_1),(A_2,M_2))$ is the composite
\[
((A_1,M_1),(A_2,M_2))\to (A_1\otimes A_2, M_1\otimes M_2)\to (B_1\otimes B_2,\psi_!(M_1\otimes M_2))
\]
where $\psi_!(M_1\otimes M_2)\simeq M_1\otimes M_2$ in $\C$. One can further check that the $B_1\otimes B_2$-coaction on $\psi_!(M_1\otimes M_2)$ must be given by the composite:
\begin{align*}
M_1\otimes M_2\to &M_1\otimes A_1\otimes M_2\otimes A_2\\\simeq & M_1\otimes M_2\otimes A_1\otimes A_2\overset{id\otimes\psi}\longrightarrow  M_1\otimes M_2\otimes B_1\otimes B_2 
\end{align*}
as in the classical case.
\end{rem}

\begin{thm}\label{thm:bialgebramonoidalstructure}
Let $\C$ be an $\E_n$-monoidal $\infty$-category for $n>1$ and let $H$ be a $(j,k)$-bialgebra in $\C$, with  $j,k\geq 1$ and $j+k\leq n$. Then $LCoMod_H(\C)$ admits the structure of an $\E_k$-monoidal $\infty$-category. Moreover, the forgetful functor $U\colon LCoMod_H(\C)\to \C$ is $\E_k$-monoidal.
\end{thm}

\begin{proof}

Since $H$ is by definition an $\E_k$-algebra in $CoAlg_{\E_j}(\C)$ there is a morphism of $\infty$-operads $H^\otimes\colon\E_k^\otimes\to CoAlg_{\E_j}(\C)^\otimes$. Pulling back the cocartesian fibration $f\colon LCoMod_{\E_j}(\C)^\otimes\to CoAlg_{\E_j}(\C)^\otimes$ of Proposition \ref{prop:cocartlmodfibration} along $H^\otimes$ gives a cocartesian fibration $p\colon LCoMod_H(\C)^\otimes\to \E_k^\otimes$. 

To see that the forgetful functor $U^\otimes\colon LCoMod_H(\C)^\otimes\hookrightarrow LCoMod_{\E_j}(\C)^\otimes\to\C^\otimes$  is $\E_k$-monoidal, first notice that, by \cite{htt}*{2.4.1.12, 2.4.2.8}, $LCoMod_H(\C)^\otimes\to LCoMod_{\E_j}(\C)^\otimes$ preserves and reflects $p$-cocartesian morphisms so therefore determines an $\E_k$-monoidal functor. By Lemma \ref{lemma:monoidalstructureonLMod}, the functor $LCoMod_{\E_j}(\C)\to \C$ is $\E_k$-monoidal.

\end{proof}
\noindent By taking opposite categories we obtain the following dual result:

\begin{cor}\label{cor:monoidalstructureonbialgebramodules}
Let $\C$ be an $\E_n$-monoidal $\infty$-category for $n\geq 2$ and let $H$ be a $(j,k)$-bialgebra in $\C$, with  $j,k\geq 1$ and $j+k\leq n$. Then $LMod_H(\C)$ admits the structure of an $\E_j$-monoidal $\infty$-category. Moreover, the forgetful functor $U\colon LMod_H(\C)\to \C$ is $\E_j$-monoidal.
\end{cor}

\begin{rem}
By applying the explicit description of cocartesian morphisms in the fibration $LCoMod_{\E_j}(\C)\to CoAlg_{\E_k}(\C)$ given in Remark \ref{expose:cocartmorphsincomod}, and the fact that $LCoMod_H(\C)\to \C$ is $\E_k$-monoidal, one can deduce that the tensor product of two $H$-comodules $M$ and $N$, using the $\E_k$-monoidal structure of Theorem \ref{thm:bialgebramonoidalstructure}, has underlying $\C$-object $M\otimes N$ with $H$-coaction map the composite \[ M\otimes N\to M\otimes H\otimes N \otimes H\simeq M\otimes N\otimes H\otimes H\to M\otimes N\otimes H \] where the final morphism is the multiplication of $H$. Similarly, if $M$ and $N$ are two $H$-modules then their tensor product is $M\otimes N$ and the $H$-action map is \[M\otimes N\otimes H\to M\otimes N\otimes H\otimes H\simeq M\otimes H\otimes N\otimes H\to M\otimes N\] where the first map is the diagonal map of $H$. Note that the requirement that $\C$ be at least $\E_2$-monoidal arises in two ways here: first, it is necessary for defining a bialgebra $H$, and second, it is required to ``twist'' the tensor product of the middle term in both of the above composites.
\end{rem}

\begin{rem}
	There is an alternative proof of Theorem \ref{thm:bialgebramonoidalstructure} in the case that $\C=\Sp$ sketched to the author by Maxime Ramzi which goes, roughly, as follows: notice that the functor $Perf\colon Alg_{\E_k}(\Sp)\to Cat_\infty^{perf}$ (where $Cat_\infty^{perf}$ is the $\infty$-category of idempotent-complete stable $\infty$-categories), which takes an algebra $A$ to  its category of perfect (equivalently compact) right modules $RMod_A^{perf}\!\!(\Sp)$, is symmetric monoidal. In particular, given a $(j,k)$-bialgebra $H\in Coalg_{\E_j}(Alg_{\E_k}(\Sp))$, the $\infty$-category $Perf(H)\simeq RMod_H^{perf}\!\!(\Sp)$ is an $\E_j$-coalgebra in $Cat_\infty^{perf}$. Now, consider the symmetric monoidal functor  $Fun^{ex}(-,\Sp)\colon Cat_\infty^{perf}\to Pr^{L,st}$ which takes a stable $\infty$-category to its $\infty$-category of exact functors to $\Sp$. Then the proof of \cite{ha}*{7.2.4.3} implies that $Fun^{ex}(RMod_H^{perf}\!\!(\Sp),\Sp)\simeq LMod_H(\Sp)$ and therefore the $\E_j$-coalgebra structure of $RMod_H^{perf}\!\!(\Sp)$ induces an $\E_j$-monoidal structure on $LMod_H(\Sp)$. A dual argument gives the result for $LCoMod_H(\Sp)$. It is likely that the above argument can be generalized to categories other than $\Sp$, and that it produces an equivalent monoidal structure to those described in Theorem \ref{thm:bialgebramonoidalstructure} and Corollary \ref{cor:monoidalstructureonbialgebramodules}, but we will not pursue that here. As always, any mistakes or errors in the above argument are the author's.
\end{rem}

\section{Some Examples of Coalgebraic Structure}

\subsection{Comonads and Descent Data}\label{sec:descendata}

We now describe how to obtain coalgebras from comonads. This is essentially an application of an Eilenberg-Watts type theorem, where we recognize comonads as coalgebras in endofunctor categories and produce coalgebras in the source category by evaluating at the generating object. This procedure allows one to recognize $\infty$-categories of descent data as categories of comodules over a coalgebra, as in \cite{hess}.

\begin{defi}
For any $\infty$-category $\C$ there is an $\E_1$-monoidal category of functors $Fun(\C,\C)$, where the monoidal structure is given by composition (cf. Remark 4.7.2.31 of \cite{ha}). A comonad on $\C$ is an object of $CoAlg(Fun(\C,\C))$.
\end{defi}

\begin{thm}\label{thm:monoidalEW}[Eilenberg-Watts]
Let $R$ be an $\E_2$-algebra in $\Sp$, and let $A$ be an $\E_1$-algebra of $LMod_R(\Sp)$.  Then there is an equivalence of monoidal categories  \[{}_ABiMod_A(LMod_R(\Sp))\overset{\sim}\to Fun^L(LMod_A(LMod_R(\Sp)),LMod_A(LMod_R(\Sp)))\] between the category of $(A,A)$-bimodules in $LMod_R(\Sp)$, and the category of colimit preserving endofunctors of left $A$-modules in $LMod_R(\Sp)$. On objects the equivalence is given by $M\mapsto M\otimes_A -$.  Moreover the respective actions of each on $LMod_A(LMod_R(\Sp))$ are equivalent. 
\end{thm}

\begin{proof}
First note that, by \cite[7.1.3.1]{ha}, it suffices to take $R= \mathbb{S}$. Now let $\mathfrak{Alg}$ be the $(\infty,2)$-category of \cite[Definition 4.40]{haugsengmorita} whose objects are $\E_1$-ring spectra and in which the $\infty$-category of morphisms from $A$ to $B$ is ${}_ABiMod_B(\Sp)$. In \cite{haugsengmorita} $\mathfrak{Alg}$ is constructed as a complete Segal object, but we pass to its equivalent representation as a $\Cat$-enriched $\infty$-category in the sense of \cite{gepnerhaugsengenriched} (i.e.~an $(\infty,2)$-category) via \cite[Corollary 7.21]{haugsengrectification} . Now from \cite[Equation 136]{hinichyoneda}, by taking $\mathcal{M}=Cat_\infty$ and $A=\mathfrak{Alg}$, we have a $Cat_\infty\times Cat_\infty$-enriched ``hom'' functor $\mathfrak{Y}\colon \mathfrak{Alg}\times\mathfrak{Alg}^{op}\to \mathfrak{Cat}$, where $\mathfrak{Cat}$ is the $(\infty,2)$-category of $(\infty,1)$-categories (as constructed by taking $\mathcal{V}$ to be $Cat_\infty$ in \cite[Definition 7.4.10]{gepnerhaugsengenriched}). Note that the hom functor above, equivalently the $(\infty,2)$-categorical Yoneda embedding, does not depend on the model of $(\infty,2)$-categories as a result of \cite{macphersonYoneda}.  By construction, the functor $\mathfrak{Y}$ takes a pair $(A,B)$ to the $\infty$-category ${}_BBiMod_A(\Sp)$ and takes a pair of modules $(M,N)\in {}_ABiMod_C(\Sp)\times{}_DBiMod_B(\Sp)$ to the functor ${}_BBiMod_A(\Sp)\to {}_DBiMod_C(\Sp)$ given by $P\mapsto N\otimes_B P\otimes_A M$.

Now let $\star_{\S}$ be the full sub-$(\infty,2)$-category of $\mathfrak{Alg}$ spanned by $\S$, i.e.~a $\Cat$-enriched $\infty$-category with one object whose category of endomorphisms is $\Sp$. By including $\star_{\S}^{op}$ into $\mathfrak{Alg}^{op}$, we obtain a functor $\mathfrak{Alg}\times\star_{\S}^{op}\to \mathfrak{Cat}$. By adjunction we have a $\Cat$-enriched functor $\overline{\mathfrak{Y}}\colon \mathfrak{Alg}\to \mathfrak{Fun}(\star_\S^{op},\mathfrak{Cat})$ from $\mathfrak{Alg}$ to the $\infty$-category of $\Cat$-enriched presheaves on $\star_{\S}$. Note that the fact that $\mathfrak{Fun}(\star_\S^{op},\mathfrak{Cat})$ is itself $\Cat$-enriched follows from \cite[Example 7.4.11]{gepnerhaugsengenriched}. By unwinding the definitions one sees that $\overline{\mathfrak{Y}}(A)$ is the functor that takes the unique object of $\star_\S$ to ${}_ABiMod_{\S}(\Sp)\simeq LMod_A(\Sp)$ and which takes each $X\in\Sp$ to $-\otimes_\S X$. Being a $\infty$-category enriched functor, this gives a  functor of mapping $\infty$-categories \[End_{\mathfrak{Alg}}(A)\simeq {}_ABiMod_{A}(\Sp)\to Fun(LMod_A(\Sp),LMod_A(\Sp))\simeq End_{\mathfrak{Cat}_\infty}(\mathfrak{I}^\ast\overline{\mathfrak{Y}}(A))\] 

 Since $LMod_A(\Sp)$ is presentable and tensoring preserves colimits (and the inclusion $\mathfrak{I}\colon \mathfrak{Pr}^L\to \mathfrak{Cat}_\infty$ is full on mapping $\infty$-categories, i.e.~$\mathfrak{Pr}^L(\C,\dD)\to \mathfrak{Cat}_\infty(\C,\dD)$ is the inclusion of a full sub-$\infty$-category), we have a lift $\mathfrak{I}^\ast\overline{\mathfrak{Y}}\colon \mathfrak{Alg}\to \mathfrak{Fun}(\star_\S^{op},\mathfrak{Pr}^L)$ where $\mathfrak{Pr}^L$ is the $(\infty,2)$-category of presentable $\infty$-categories and left adjoints between them (the existence of such an $(\infty,2)$-category is implied by first noticing that the $\infty$-category $Pr^L$ is enriched over itself by \cite[7.4.10]{gepnerhaugsengenriched}, and then applying \cite[4.3.9]{gepnerhaugsengenriched} to the lax monoidal forgetful functor of $\infty$-categories $Pr^L\to \Cat$). 

It follows from \cite[Definition 2.2.17]{gepnerhaugsengenriched} that the functor of mapping $\infty$-categories induced by $\mathfrak{I}^\ast\overline{\mathfrak{Y}}$, namely \[End_{\mathfrak{Alg}}(A)\simeq {}_ABiMod_{A}(\Sp)\to Fun^L(LMod_A(\Sp),LMod_A(\Sp))\simeq End_{\mathfrak{Pr}^L}(\mathfrak{I}^\ast\overline{\mathfrak{Y}}(A))\] is a functor of $\mathbb{E}_1$-monoidal $\infty$-categories (where the monoidal structure is ``composition'').  Additionally, it follows that their respective actions on $LMod_A(\Sp)$ are compatible with this monoidal equivalence (again by ``composition'').  By \cite[7.1.2.4]{ha} the functor ${}_ABiMod_A(\Sp)\to Fun^L(LMod_A(\Sp),LMod_A(\Sp))$ is an equivalence, completing the proof.
\end{proof}

The following is an $\infty$-categorical version of \cite[25.4]{brz} and is inspired by \cite{hess}.

\begin{thm}\label{thm:descentdatacoring}
Let $R$ be an $\E_2$-ring spectrum and $\phi\colon A\to B$ be a morphism of $\E_1$-$R$-algebras. Then the $\infty$-category of descent data for $\phi$ is equivalent to the $\infty$-category of comodules in $LMod_B(\Sp)$ for the descent coring $B\otimes_A B$.
\end{thm}

\begin{proof}
By \cite[4.7.5.2 (3)]{ha}, the $\infty$-category of descent data for a morphism of $\E_1$-ring spectra $\phi\colon A\to B$ is equivalent to the $\infty$-category of coalgebras over a comonad on $LMod_B(\Sp)$, and this comonad is given by tensoring over $B$ with the $(B,B)$-bimodule $B\otimes_A B$. By Theorem \ref{thm:monoidalEW}, this comonad determines, and is determined by, a coalgebra in ${}_BBiMod_B(\Sp)$, namely $B\otimes_A B$. Moreover, again by Theorem \ref{thm:monoidalEW}, the category of coalgebras in $LMod_B(\Sp)$ over the comonad $-\otimes_B B\otimes_A B$ is equivalent to the category of $B\otimes_A B$-comodules in $LMod_B(\Sp)$. 
\end{proof}

\begin{rem}
Theorem \ref{thm:descentdatacoring} can be thought of as a weak form of ``Tannakian realization'' for $LMod_A(\Sp)$ in the sense that, if $\phi$ is an effective descent morphism, then  $LMod_A(\Sp)$ is equivalent to the category of descent data for $\phi$ and is therefore equivalent to a category of comodules over a coalgebra. In the special case of a Galois or Hopf-Galois extension of ring spectra (cf.~\cite{rog}), we have respective equivalences $B\otimes_A B\simeq B\otimes G^{\vee}$ and $B\otimes_A B\simeq B\otimes H$ where $G$ is a stably dualizable group (again in the sense of \cite{rog}) and $H$ is a spectral bialgebra. In this case we can take the Tannakian interpretation further and think of $LMod_A(\Sp)$ as the category of representations of an honest (non-commutative, spectral) group scheme. Moreover, descent along the morphism $\phi$ can respectively be reinterpreted as taking fixed points with respect to the $G$-action, or primitives with respect to the $H$-coaction. These give, respectively, a homotopy fixed points spectral sequence and an Adams spectral sequence. 
\end{rem}

\begin{rem}
	The above theorem induces a number of questions that deserve a more detailed analysis. For instance, what \textit{precisely} is the coalgebra structure of $B\otimes_AB$ in ${}_BBiMod_B(\Sp)$? In \cite{hess} it is defined to be the composite \[B\otimes_AB\overset{\sim}\to B\otimes_A A\otimes_A B\overset{1\otimes\phi\otimes 1}\to B\otimes_A B\otimes_A B\overset{\sim}\to  B\otimes_AB\otimes_B B\otimes_AB  \] and one checks that this determines a comonoid structure and that descent data are comodules for that comonoid. In our $\infty$-categorical setting however, the analysis goes in the reverse direction. Namely, one proves that there \textit{is} a coalgebra controlling descent whose underlying spectrum is $B\otimes_AB$, and then one needs to investigate precisely which morphism in $({}_BBiMod_B(\Sp)^{op})^\otimes$ (determined by some cocartesian lift of a map in $Ass^\otimes$) gives the coalgebra structure (an analysis which is similar in spirit to our Remark \ref{expose:cocartmorphsincomod}).   
\end{rem}

\subsection{Thom Spectra}\label{sec:thomobjects}

In what follows, we primarily use the language and constructions of \cite{abg,  abghr,  acb}. Recall that $Pic\colon Alg_{\Op}(Pr^L)\to Alg_{\Op}^{gp}(\T)$ is the functor that takes an $\Op$-monoidal presentable $\infty$-category to its $\Op$-monoidal and grouplike Picard space of invertible objects and equivalences between them (as in \cite{abg}). In particular we use the standard notation of parameterized homotopy theory (as in e.g.~\cite{abg,abghr}), and we will write $\C^X$ to denote the $\infty$-category of functors $Fun(X,\C)$. 

 If $R$ is an $\E_{n}$-ring spectrum, we will write $Pic(R)$, instead of $Pic(LMod_R(\Sp))$, for its $(n-1)$-fold loop space of invertible objects. Note that there is always an inclusion $Pic(R)\hookrightarrow LMod_R(\Sp)$. 

All of the results of this section would hold if we replaced $\Sp$ by some other presentable stable $\infty$-category (and most would hold even if we replaced $\Sp$ by some other presentable but not necsesarily stable $\infty$-category), but we work with $\Sp$ for the reason that the examples of interest are all contained therein. To simplify notation, we will not write the category in which our algebras and modules live (since it will always be $\Sp$), e.g.~we will simply write $Alg_{\E_k}$ and $LMod_R$ rather than $Alg_{\E_k}(\Sp)$ or $LMod_R(\Sp)$.

\begin{defi}[\cite{abghr}*{Definition 1.4}]\label{def:thomfunctor}
	For $R$ an $\E_n$-ring spectrum we define the Thom spectrum functor $\Th(-)\colon\T_{/Pic(R)}\to LMod_R$ to be $colim_X(i\circ -)$, where $i\colon Pic(R)\hookrightarrow LMod_R$ is the inclusion. 
\end{defi}

The following is originally due to Lewis, but we use a modern reference:

\begin{thm}[\cite{acb}*{Theorem 2.8}]\label{thm:acbmonoidal}
	If $f\colon X\to Pic(R)$ is a map of $\E_k$-monoidal $\infty$-groupoids (hence the composite $i\circ f$ is an $\E_k$-monoidal functor) for $0\leq k\leq n-1$, then $\Th(f)$ is an $\E_k$-algebra in $LMod_R$. 
\end{thm}

\begin{lem}\label{lem:trivialthom}
	Let $R$ be an $\E_{n}$-ring spectrum for $n>0$. If $R_X\colon X\to Pic(R)$ is any map that factors through a contractible simplicial set, hence equivalent to the constant functor valued in $R$, then $\Th(R_X)\simeq R\otimes_R\otimes R\otimes  \Sigma_+^\infty X\simeq R\otimes \Sigma^\infty_+X$.
\end{lem}

\begin{proof}

	This is \cite{abghr}*{Proposition 2.8}.
\end{proof}

From Lemma \ref{lem:trivialthom} it follows that the Thom functor applied to the constant morphism always produces a coalgebra object:

\begin{prop}\label{thm:trivialthomiscoalg}
	Let $R$ be an $\E_{n}$-algebra in $\Sp$ for $n>0$ and let $R_X\colon X\to Pic(R)$ be a constant map valued in $R$. Then $\Th(R_X)$ is an $\E_{n-1}$-coalgebra in $LMod_R$ via the diagonal map of $X$. 
\end{prop}

\begin{proof}
	Both $\Sigma^\infty_+\colon \T\to\Sp$ and $-\otimes R\colon \Sp\to LMod_R$ are $\E_n$-monoidal so an application of Proposition \ref{prop:functorofbialgebra} completes the proof.
\end{proof}

\begin{rem}
	Note that while $X$ is an $\E_\infty$-coalgebra in $\T$, it cannot have the same degree of cocommutativity in $LMod_R$, as the latter is only an $\E_{n-1}$-monoidal $\infty$-category.
\end{rem}

\subsection{The Thom Diagonal and Thom Isomorphism}\label{sec:thomarecomods}

We now wish to show that for any $\E_{n}$-ring spectrum $R$, with $n>1$, and any morphism $f\colon X\to Pic(R)$, the $R$-module $\Th(f)$ is a comodule over the $\E_{n-1}$-coalgebra $\Th(R_X)$. Later we will show (cf.~Theorem \ref{thm:thomdiagonal}) that this coaction is equivalent to the ``classical'' Thom diagonal as described, for instance, in \cite{abghr} and that it participates in the Thom isomorphism in the usual way (cf.~Theorem \ref{thm:thomisomodules}). Moreover, when $f$ is a map of $\E_k$-monoidal $\infty$-groupoids then $\Th(f)$ is an $\E_k$-$R$-algebra and this coaction can be used to imbue $\Th(f)$ with a kind of non-commutative $R\otimes\Sigma^\infty_+ X$-torsor structure over $R$.

\begin{rem}
Because we are assuming $R$ is at least an $\E_2$-ring spectrum, we will be somewhat indelicate about the difference between left modules and right modules. 
\end{rem}

\begin{lem}\label{lem:oplax}
	Let $\C$ be an $\mathbb{E}_k$-monoidal complete and cocomplete $\infty$-category and $X$ a Kan complex with terminal morphism $p\colon X\to \ast$. Then the colimit functor $p_!\colon \C^X\to \C$ is oplax $\E_k$-monoidal and the limit functor $p_\ast\colon \C^X\to \C$ is lax $\E_k$-monoidal, in the sense of \cite{haugsengmates}, where $\C^X$ is equipped with the pointwise monoidal structure.
\end{lem}

\begin{proof}
	The colimit and limit functors are left and right adjoint, respectively, to the diagonal functor $p^\ast\colon \C\to \C^X$. The diagonal functor is (strongly) $\mathbb{E}_k$-monoidal by \cite{abg}*{Theorem 6.4} and therefore both lax monoidal and oplax monoidal. The result then follows from \cite{haugsengmates}*{4.5, 4.6}.  
\end{proof}

\begin{rem}
	Notice that, with respect to the pointwise monoidal structure on $LMod_R^X$, $R_X$ is the monoidal unit. As such, it is a coalgebra and every $f\in LMod_R^X$ is a right (and left) comodule over it.
\end{rem}

\begin{cor}\label{thm:thomdiagonal}
	For a functor of $\infty$-groupoids $f\colon X\to Pic(R)$, $\Th(f)$ is a $\Th(R_X)\simeq R\otimes \Sigma^\infty X$-comodule in $LMod_R(\Sp)$. In particular, there is a coassociative coaction $\Delta_f\colon\Th(f)\to R\otimes\Sigma^\infty_+X\otimes_R\Th(f)\simeq \Sigma^\infty_+ X\otimes \Th(f) $ in $LMod_R$.
\end{cor}

\begin{thm}\label{thm:thomisomodules}
Let $f\colon X\to Pic(R)$ be a functor of $\infty$-groupoids and suppose that $\Th(f)$ has an $L$-type $E$-orientation $\theta\colon \Th(f)\to E$ for some $E\in Alg_{\E_k}(LMod_R)$. Then the Thom isomorphism of \cite{abghr}*{Corollary 2.26} decomposes as the composite \[E\otimes_R \Th(f)\overset{1_E\otimes \Delta_f}\to E\otimes_R \Th(f)\otimes \Sigma^\infty_+X \overset{1_E\otimes \theta\otimes 1_X}\to E\otimes_R E\otimes X\overset{\mu_E\otimes 1_X}\to E\otimes X \] in which the first map is the coaction map of Corollary \ref{thm:thomdiagonal} tensored with the identity on $E$, the second map uses the orientation on the middle component, and the final map is the multiplication of $E$.
\end{thm}

\begin{proof}
We begin by describing the Thom isomorphism of \cite{abghr}*{2.26} in detail. For the extent of this proof we will denote the terminal morphism for $X$ in $\T$ by $p\colon X\to \ast$, and write $p_!$ for the colimit functor $LMod_R^X$ and $p^\ast$ for the diagonal $LMod_R\to LMod_R^X$.  By assumption we have a morphism $\theta\colon \Th(f)\to E$ in $\C$ which induces an equivalence on adjoints when restricted to points of $X$. We obtain by Proposition \ref{prop:orientationequivalence} an orientation in the sense of \cite{abghr}, i.e.~a morphism $\tilde{\theta}\colon \Th(f)\otimes E\to E$ in $LMod_E(\C)$ which decomposes as $p_!(f\otimes E)\simeq\Th(f)\otimes E\overset{\theta\otimes E}\to E\otimes E\overset{\mu_E}\to E$ and whose adjoint in $LMod_E^X$ is an equivalence. From the proof of \cite{htt}*{5.2.2.8} we have that this adjoint of $\tilde{\theta}$ is precisely the composite
\[
f\otimes E\to (p^\ast p_! f)\otimes E\simeq p^\ast (\Th(f)\otimes E)\to p^\ast (E\otimes E)\to p^\ast E
\] 
in which the the first map is the unit of the $p_!\dashv p^\ast$-adjunction (followed by a straightforward equivalence), the second is $p^\ast (\theta\otimes E)$ and the last is $p^\ast \mu_E$. The Thom isomorphism of \cite{abghr}*{Corollary 2.26} is defined to be the colimit of the above composite. Note that applying $p_!$ to the second and third morphisms above, which are morphisms of constant functors, gives the composite $X\otimes\Th(f)\otimes E\to X\otimes E\otimes E\to X\otimes E$ as desired. It remains to show that the first map, namely $p_!$ of the unit $f\to p^\ast p_!f$ (tensored with $E$), is equivalent to the oplax structure map $\Delta_f$ whose existence is guaranteed by Lemma \ref{lem:oplax}.

From \cite{HLN-orthofibrations}*{7.7} (or by the proof of \cite{haugsengmates}*{2.2}) we have an explicit formula for the result of applying an oplax monoidal functor to a coaction. Thus we can determine a formula for $\Delta_f$, which arises as the colimit of the (identity) coaction $f\overset{\sim}\to f\otimes R_X$ in $LMod_R^X$. In general, for a fixed pair of functors $F,G\in LMod_R^X$, the preceding two references imply that the oplax structure morphism is  \[p_!(F \otimes G)\to p_!(p^\ast p_!(F) \otimes p^\ast p_!(G)) \overset{\sim}\to p_! p^\ast(p_! F \otimes p_! G) \to p_! F \otimes p_! G\] where the first morphism is the tensor of two units of the adjunction, the second morphism is the lax monoidal structure map of $p^\ast$ and the final map is the counit. Note that because $p^\ast$ is in fact strong monoidal the second morphism above is an equivalence. This morphism can be decomposed further as \[p_!(F \otimes G)\to p_!(p^\ast p_!(F) \otimes G)\to p_!(p^\ast p_!(F) \otimes p^\ast p_!(G))  \overset{\sim}\to p_! p^\ast(p_! F \otimes p_! G) \to p_! F \otimes p_! G \] by applying the unit first to $F$ and then to $G$.  One can check that the composite \[p_!(p^\ast p_!(F) \otimes G)\to p_!(p^\ast p_!(F) \otimes p^\ast p_!(G))  \overset{\sim}\to p_! p^\ast(p_! F \otimes p_! G) \to p_! F \otimes p_! G\] is an equivalence because colimits in $LMod_R$ commute with tensor product, (it is in fact a decomposition of the so-called ``projection formula'' of parameterized homotopy theory). As a result, by taking $G=R_X$ and $F=f$, and using that $R_X$ is the monoidal unit of $LMod_R^X$, we have that $\Delta_f$ is the following composite: 
\[
\Th(f)\simeq p_!(f\otimes R_X)\simeq p_!f\to p_!(p^\ast p_! f)\simeq p_!(p^\ast p_! f\otimes R_X)\simeq p_!f\otimes_R p_!R_X\simeq \Th(f)\otimes \Sigma^\infty_+ X
\]
and is therefore equivalent to the colimit of the application of the unit $id\to p^\ast p_!$ to $f$, which is the first map of the first displayed composite of this proof.
\end{proof}

\begin{cor}
	If $f\colon X\to Pic(R)$ is $\E_k$-monoidal for $1\leq k\leq n-1$ then the following composite is an equivalence: \[ \Th(f)\otimes_R\Th(f)\overset{1_{\Th(f)}\otimes \Delta_f}\to \Th(f)\otimes_R\Th(f)\otimes\Sigma^\infty_+X\overset{\mu_{\Th(f)}\otimes 1_X}\to \Th(f)\otimes \Sigma^\infty_+X \]
\end{cor}

\begin{proof}
	This follows from the fact that $\Th(f)$ is always oriented with respect to itself via the identity map, by \cite{acb}*{Corollary 3.17}.
\end{proof}

\appendix

\section{Three equivalent notions of orientation}\label{appendix: orientations}

In the following appendix, we review the notion of an \textit{orientation} for a Thom spectrum with respect to a ring spectrum. In the preceding sections, we use the results of this appendix sparingly (indeed  exactly once, in Theorem \ref{thm:thomisomodules}). However they still may be useful to record, if only as an aid to those new to the concept. We continue with the same notation and terminology from Section \ref{sec:thomobjects}. 

The below is an $\infty$-categorical reformulation of the notion of $E$-orientation described by Lewis in \cite{lewis}.

\begin{defi}\label{def:Lewisorientation}
Let $f\colon X\to Pic(R)$ be a morphism of $\infty$-groupoids for $R$ some  $\E_n$-ring spectrum with $n\geq 2$. Let $E$ be an $\E_k$-$R$-algebra for $1\leq k\leq n-1$. Then we say that an $L$-type $E$-orientation of $\Th(f)$ is the data of an $R$-module map $u\colon \Th(f)\to E$ satisfying the following property for all points $x\in X$:

\begin{itemize}
\item Let $\Th(x)$ be the Thom spectrum associated to the inclusion $x\hookrightarrow X\to Pic(R)$ and let $u_x\colon \Th(x)\to E$ by the composition of the induced map $\Th(x)\to \Th(f)$ with $u\colon \Th(f)\to E$. Then the image of $u_x$ under the equivalence of mapping spaces induced by the extension of scalars adjunction from $R$ to $E$, $LMod_R(\Th(x),E)\overset{\sim}\to LMod_E(E\otimes \Th(x),E)$, i.e.~the composition $E\simeq \Th(x)\otimes E\overset{u_x\otimes E}\to E\otimes E\overset{\mu_E}\to E$, is an equivalence of $E$-modules. 
\end{itemize}
\end{defi}

The below is the definition of $E$-orientation given in \cite{abghr}.

\begin{defi}\label{def:ABGHRorientation}
Let $E$ and $R$ be as in Definition \ref{def:Lewisorientation}. Let $\Th(g)$ be the Thom spectrum associated to a map of $\infty$-groupoids $g\colon X\to Pic(E)$. Then an $A$-type $E$-orientation of $\Th(g)$ is a morphism of $E$-modules $u\colon \Th(g)\to E$ with the following property:
\begin{itemize}
\item Let $p\colon X\to \ast$ be the terminal map, which induces an adjunction $p_!\colon LMod_E^X\rightleftarrows LMod_E\colon p^\ast$, and note that by definition $\Th(g)\simeq p_!(g)$. Then the image of  $u$ under the equivalence $LMod_E(p_!(f),E)\overset{\sim}\to LMod_E^X(f,p^\ast(E))$ is an equivalence. 
\end{itemize}
\end{defi}

\begin{lem}\label{lemma:sliceequivalence}
Let $\C$ be be an $\infty$-category with $E\in \C$ and $X$ a simplicial set. Then there is an equivalence of $\infty$-categories $(\C_{/E})^X\simeq (\C^X)_{/E_X}$ where $E_X$ is the constant functor $X\to \C$ taking every point of $X$ to $E$.
\end{lem}

\begin{proof}
From \cite{rvelements}*{4.2.1} the slice $\infty$-category $\C_{/E}$ is defined by the following pullback square \[
\xymatrix{
\C_{/E}\ar[r]\ar[d] & \C^{\Delta^1}\ar[d]^{(p_0,p_1)}\\
\C\times \Delta^0\ar[r]^{1_\C\times E} & C\times C
}
\] in which the right hand vertical map is induced by the inclusion of the endpoints $\Delta^0\times\Delta^0\hookrightarrow \Delta^1$. Note that, because it is obtained by cotensoring with the above monomorphism, the projection $(p_0,p_1)\colon\C^{\Delta^1}\to\C\times \C$ is an isofibration by definition (cf.~\cite{rvelements}*{Definition 1.2.1 (ii)}) and therefore pulling back along it is a cosmological limit by \cite{rvelements}*{1.2.1(i)} and the ensuing discussion. From \cite{rvelements}*{1.3.4(iii)} we have that the functor $(-)^X$ preserves the displayed pullback square (because $(-)^X$ is a cosmological functor which by \cite{rvelements}*{1.3.1} preserves cosmological limits). As a result we have another pullback square
\[
\xymatrix{
(C_{/E})^X\ar[r]\ar[d] & (\C^X)^{\Delta^1}\ar[d] \\
\C^X\times\Delta^0\ar[r] & \C^X\times\C^X
}
\]

(in which we've already applied the obvious equivalences to the bottom left, bottom right, and upper right corners). Note that the bottom horizontal functor in the above pullback diagram is the product of the identity functor $1_{\C^X}$ and the functor which picks out the constant functor $E_X\colon X\to \C$. Therefore the lower right cospan of the above square is the one whose pullback is $(C^X)_{/E_X}$, giving the desired equivalence.
\end{proof}

\begin{rem}
Note that the ``slice $\infty$-category" of \cite{htt} is not defined in the same way as that of \cite{rvelements}. However by \cite{rvelements}*{D.6.6} the two notions are equivalent $\infty$-categories and so Lemma \ref{lemma:sliceequivalence} is applicable to the slice $\infty$-categories used in this paper. 
\end{rem}

\begin{rem}\label{remark:corefactorization}
Recall that, given a natural transformation $\eta\colon F\Rightarrow G$ of functors $F,G\in \C^X$, $\eta$ is an equivalence of functors if and only if the restrictions $\eta_x\colon F_x\Rightarrow G_x$ are equivalences for all $x\in X$. As a result the equivalence of Lemma \ref{lemma:sliceequivalence} identifies functors $X\to C_{/E}$ that factor through the full subcategory of $\C_{/E}$ spanned by equivalences to $E$ with natural transformations to $E_X$ which are equivalences (note that these are \textit{not} the cores of $\C_{/E}$ or $C^X$ but rather full subcategories spanned by certain classes of \textit{objects}). 
\end{rem}

\begin{prop}\label{prop:orientationequivalence}
Let $f$, $R$ and $E$ be as in Definition \ref{def:Lewisorientation}. Then a map of $R$-modules $\theta\colon \Th(f)\to E$ is an $L$-type $E$orientation if and only its adjoint $\tilde{\theta}\colon \Th(f)\otimes E\to E$ is an $A$-type $E$-orientation. 
\end{prop}

\begin{proof}
First suppose that $\theta\colon \Th(f)\to E$ is an $L$-type orientation for a Thom spectrum associated to a map $f\colon X\to Pic(R)$. Then by using the relevant adjunctions we obtain equivalences of mapping spaces $LMod_R(p_!(f),E)\simeq LMod_R^X(f, p^\ast(E))\simeq LMod_R^X(f, E_X)$ so that $\theta$ determines an object of $(LMod_R^X)_{/E_X}$ (which in a slight abuse of notation we will also call $\theta$), namely the natural transformation which on components if given by the composite $R\to \Th(f)\overset{\theta}\to E$, where the first map is the unique map to the colimit. Therefore by Lemma \ref{lemma:sliceequivalence} we have a functor $f_{\theta}\colon X\to (LMod_R)_{/E}$ which takes each point of $x$ to the component wise map given above. By applying \cite{htt}*{5.2.5.1} to the extension/restriction of scalars adjunction between $LMod_R$ and $LMod_E$, this determines a functor $\tilde{f}_{\theta}\colon X\to (LMod_E)_{/E}$ (or equivalently by applying the free $E$-module functor and then composing with the multiplication of $E$). Note that because colimits in overcategories are computed in the underlying category (cf.~\cite{htt}*{1.2.13.8}), the colimit of $\tilde{f}_{\theta}$ in $(LMod_E)_{/E}$ is $\tilde{\theta}\colon \Th(f)\otimes E\to E$, the adjoint of the Lewis orientation $\theta\colon \Th(f)\to E$. Similarly, for each $x\in X$ the induced functor $x\hookrightarrow X\to EMod_{/E}$ is $\tilde{\theta}_x\colon \Th(x)\otimes E\to E$, the adjoint of the restriction of the Lewis orientation $\theta_x\colon \Th(x)\to E$. By assumption, each morphism $\tilde{\theta}_x$ is an equivalence. Therefore by Lemma \ref{lemma:sliceequivalence} and Remark \ref{remark:corefactorization} this induces an equivalence of functors $\tilde{f}_\theta\simeq E_X\in LMod_E^X$. Therefore $\tilde{\theta}$ is an $A$-type orientation of $\Th(f)\otimes E$. The reverse implication follows from the above argument in reverse. 
\end{proof}

\begin{cor}
Let $f\colon X\to Pic(R)$ be a morphism of $\infty$-groupoids with associated Thom spectrum $\Th(f)$. Then a map of $R$-modules $\theta\colon \Th(f)\to E$ is an $L$-type $E$-orientation if and only if its adjoint $\tilde{\theta}\colon \Th(f)\otimes E\to E$ is an $A$-type $E$-orientation if and only if the composite $X\overset{f}\to Pic(R)\overset{-\otimes E}\to Pic(E)$ is null. 
\end{cor}

\begin{proof}
This follows from \cite{abghr}*{2.24} which identifies $A$-type $E$-orientations with such nullhomotopies.
\end{proof}

\addcontentsline{toc}{section}{References}
\bibliography{references}

\begin{bibdiv}
\begin{biblist}

\bib{abg}{article}{
      author={Ando, Matthew},
      author={Blumberg, Andrew~J.},
      author={Gepner, David},
       title={Parametrized spectra, multiplicative {T}hom spectra and the
  twisted {U}mkehr map},
        date={2018},
        ISSN={1465-3060},
     journal={Geom. Topol.},
      volume={22},
      number={7},
       pages={3761\ndash 3825},
         url={https://doi.org/10.2140/gt.2018.22.3761},
      review={\MR{3890766}},
}

\bib{abghr}{article}{
      author={Ando, Matthew},
      author={Blumberg, Andrew~J.},
      author={Gepner, David},
      author={Hopkins, Michael~J.},
      author={Rezk, Charles},
       title={An {$\infty$}-categorical approach to {$R$}-line bundles,
  {$R$}-module {T}hom spectra, and twisted {$R$}-homology},
        date={2014},
        ISSN={1753-8416},
     journal={J. Topol.},
      volume={7},
      number={3},
       pages={869\ndash 893},
}

\bib{acb}{article}{
      author={Antol\'{\i}n-Camarena, Omar},
      author={Barthel, Tobias},
       title={A simple universal property of {T}hom ring spectra},
        date={2019},
        ISSN={1753-8416},
     journal={J. Topol.},
      volume={12},
      number={1},
       pages={56\ndash 78},
         url={https://doi.org/10.1112/topo.12084},
      review={\MR{3875978}},
}

\bib{bgn}{article}{
      author={Barwick, Clark},
      author={Glasman, Saul},
      author={Nardin, Denis},
       title={Dualizing cartesian and cocartesian fibrations},
        date={2018},
     journal={Theory Appl. Categ.},
      volume={33},
       pages={Paper No. 4, 67\ndash 94},
      review={\MR{3746613}},
}

\bib{beardsleyperoux}{misc}{
      author={Beardsley, Jonathan},
      author={Péroux, Maximilien},
       title={Koszul duality in higher topoi},
        date={2019},
        note={arXiv:1909.11724},
}

\bib{brz}{book}{
      author={Brzezinski, Tomasz},
      author={Wisbauer, Robert},
       title={Corings and comodules},
      series={London Mathematical Society Lecture Note Series},
   publisher={Cambridge University Press, Cambridge},
        date={2003},
      volume={309},
        ISBN={0-521-53931-5},
         url={https://doi.org/10.1017/CBO9780511546495},
      review={\MR{2012570}},
}

\bib{gepnerhaugsengenriched}{article}{
      author={Gepner, David},
      author={Haugseng, Rune},
       title={Enriched {$\infty$}-categories via non-symmetric
  {$\infty$}-operads},
        date={2015},
        ISSN={0001-8708},
     journal={Adv. Math.},
      volume={279},
       pages={575\ndash 716},
         url={https://doi.org/10.1016/j.aim.2015.02.007},
      review={\MR{3345192}},
}

\bib{haugsengrectification}{article}{
      author={Haugseng, Rune},
       title={Rectification of enriched {$\infty$}-categories},
        date={2015},
        ISSN={1472-2747},
     journal={Algebr. Geom. Topol.},
      volume={15},
      number={4},
       pages={1931\ndash 1982},
         url={https://doi.org/10.2140/agt.2015.15.1931},
      review={\MR{3402334}},
}

\bib{haugsengmorita}{article}{
      author={Haugseng, Rune},
       title={The higher {M}orita category of {$\Bbb{E}_n$}-algebras},
        date={2017},
        ISSN={1465-3060},
     journal={Geom. Topol.},
      volume={21},
      number={3},
       pages={1631\ndash 1730},
         url={https://doi.org/10.2140/gt.2017.21.1631},
      review={\MR{3650080}},
}

\bib{haugsengmates}{misc}{
      author={Haugseng, Rune},
       title={A fibrational mate correspondence for $\infty$-categories},
        date={2020},
        note={arXiv:2011.08808},
}

\bib{hess}{misc}{
      author={Hess, Kathryn},
       title={A general framework for homotopic descent and codescent},
        date={2010},
        note={arxiv.org/abs/1001.1556v3},
}

\bib{hinichyoneda}{article}{
      author={Hinich, Vladimir},
       title={Yoneda lemma for enriched {$\infty$}-categories},
        date={2020},
        ISSN={0001-8708},
     journal={Adv. Math.},
      volume={367},
       pages={107129, 119},
         url={https://doi.org/10.1016/j.aim.2020.107129},
      review={\MR{4080581}},
}

\bib{HLN-orthofibrations}{misc}{
      author={Hebestreit, Fabian},
      author={Linskens, Sil},
      author={Nuiten, Joost},
       title={Orthofibrations and monoidal adjunctions},
        date={2020},
        note={arXiv:2011.11042},
}

\bib{haugseng_melani_safronov}{article}{
      author={Haugseng, Rune},
      author={Melani, Valerio},
      author={Safronov, Pavel},
       title={Shifted coisotropic correspondences},
        date={2020},
     journal={Journal of the Institute of Mathematics of Jussieu},
       pages={1–65},
}

\bib{hry}{book}{
      author={Hackney, Philip},
      author={Robertson, Marcy},
      author={Yau, Donald},
       title={Infinity properads and infinity wheeled properads},
      series={Lecture Notes in Mathematics},
   publisher={Springer, Cham},
        date={2015},
      volume={2147},
}

\bib{lewis}{thesis}{
      author={Lewis, L.G.},
       title={The stable category and generalized {T}hom spectra},
        type={Ph.D. Thesis},
        date={1978},
        note={The University of Chicago},
}

\bib{htt}{book}{
      author={Lurie, Jacob},
       title={Higher topos theory},
      series={Annals of Mathematics Studies},
   publisher={Princeton University Press},
     address={Princeton, NJ},
        date={2009},
      volume={170},
}

\bib{ha}{misc}{
      author={Lurie, Jacob},
       title={Higher algebra},
        date={2014},
        note={\url{math.harvard.edu/\textasciitilde
  lurie/papers/higheralgebra.pdf}},
}

\bib{kerodon}{misc}{
      author={Lurie, Jacob},
       title={Kerodon},
        date={2021},
        note={\url{https://kerodon.net}},
}

\bib{macphersonYoneda}{article}{
      author={Macpherson, Andrew~W.},
       title={The operad that co-represents enrichment},
        date={2021},
        ISSN={1532-0073},
     journal={Homology Homotopy Appl.},
      volume={23},
      number={1},
       pages={387\ndash 401},
         url={https://doi.org/10.4310/hha.2021.v23.n1.a20},
      review={\MR{4185309}},
}

\bib{mahringthom}{article}{
      author={Mahowald, Mark},
       title={Ring spectra which are {T}hom complexes},
        date={1979},
        ISSN={0012-7094},
     journal={Duke Math. J.},
      volume={46},
      number={3},
       pages={549\ndash 559},
         url={http://projecteuclid.org/euclid.dmj/1077313574},
      review={\MR{544245}},
}

\bib{rog}{article}{
      author={Rognes, John},
       title={{Galois} extensions of structured ring spectra/{S}tably
  dualizable groups},
        date={2008},
        ISSN={0065-9266},
     journal={Mem. Amer. Math. Soc.},
      volume={192},
      number={898},
       pages={viii+137},
         url={https://doi.org/10.1090/memo/0898},
      review={\MR{2387923}},
}

\bib{froth}{thesis}{
      author={Roth, Fridolin},
       title={{Galois and Hopf-Galois} theory for associative {S}-algebras},
        type={Ph.D. Thesis},
        date={2009},
        note={Universit\"{a}t Hamburg},
}

\bib{rvelements}{misc}{
      author={Riehl, Emily},
      author={Verity, Dominic},
       title={Elements of $\infty$-category theory},
        note={Available at: https://math.jhu.edu/$\sim$eriehl/elements.pdf.
  Accessed January 20, 2021},
}

\bib{vistolistacks}{misc}{
      author={Vistoli, Angelo},
       title={Notes on {G}rothendieck topologies, fibered categories and
  descent theory},
        date={2008},
        note={homepage.sns.it/vistoli/descent.pdf},
}

\end{biblist}
\end{bibdiv}

\end{document}